
\documentstyle{amsppt}
\baselineskip18pt
\magnification=\magstep1
\pagewidth{30pc}
\pageheight{45pc}
\hyphenation{co-deter-min-ant co-deter-min-ants pa-ra-met-rised
pre-print pro-pa-gat-ing pro-pa-gate
fel-low-ship Cox-et-er dis-trib-ut-ive}
\def\leaderfill{\leaders\hbox to 1em{\hss.\hss}\hfill}
\def\A{{\Cal A}}
\def\D{{\Cal D}}
\def\H{{\Cal H}}
\def\L{{\Cal L}}
\def\Pl{{\Cal P}}

\def\afn{{\text {\bf a}}}
\def\tr{{\text {tr}}}

\def\idest{i.e.,\ }

\def\g{{\gamma}}
\def\G{{\Gamma}}
\def\d{{\delta}}

\def\e{{\varepsilon}}

\def\th{{\theta}}

\def\k{{\kappa}}
\def\l{{\lambda}}

\def\t{{\tau}}

\def\te{\widetilde t}

\def\ba{{\bold a}}

\def\bc{{\bold c}}

\def\BB{{\bold B}}
\def\b0{\text{\bf 0}}

\def\ra{{\ \longrightarrow \ }}

\def\ds{\displaystyle}

\def\supp{\text{\rm \, supp}}

\def\lan{{\langle}}
\def\ran{{\rangle}}

\def\real{{\Bbb R}}
\def\complex{{\Bbb C}}
\def\zed{{\Bbb Z}}

\def\enn{{\Bbb N}}

\def\Im{\text{\rm Im}}

\def\boxit#1{\vbox{\hrule\hbox{\vrule \kern3pt
\vbox{\kern3pt\hbox{#1}\kern3pt}\kern3pt\vrule}\hrule}}
\def\rabbit{\vbox{\hbox{\kern0pt
\vbox{\kern0pt{\hbox{---}}\kern3.5pt}}}}

\def\tableau#1{
        \hbox {
                \hskip -10pt plus0pt minus0pt
                \raise\baselineskip\hbox{
                \offinterlineskip
                \hbox{#1}}
                \hskip0.25em
        }
}

\def\tabCol#1{
\hbox{\vtop{\hrule
\halign{\strut\vrule\hskip0.5em##\hskip0.5em\hfill\vrule\cr\lower0pt
\hbox\bgroup$#1$\egroup \cr}
\hrule
} } \hskip -10.5pt plus0pt minus0pt}

\def\CR{
        $\egroup\cr
        \noalign{\hrule}
        \lower0pt\hbox\bgroup$
}



\def\blank#1#2{
\hbox to #1{\hfill \vbox to #2{\vfill}}
}


\def\strut{\vrule height10pt depth5pt width0pt}

\topmatter
\title On planar algebras arising from hypergroups
\endtitle

\author R.M. Green \endauthor
\affil 
Department of Mathematics and Statistics\\ Lancaster University\\
Lancaster LA1 4YF\\ England\\
{\it  E-mail:} r.m.green\@lancaster.ac.uk\\
\endaffil


\abstract
Let $A$ be an associative algebra with identity and with trace.  
We study the family of planar algebras on $1$-boxes that arise from
$A$ in the work of Jones, but with the added assumption that the
labels on the $1$-boxes come from a discrete hypergroup in the sense
of Sunder.  This construction equips the algebra 
$P_n^A$ with a canonical basis, $\BB_n^A$, which turns out to be a
``tabular'' basis.  We examine special cases of this construction to 
exhibit a close connection between such bases and Kazhdan--Lusztig
bases of Hecke algebras of types $A$, $B$, $H$ or $I$.
\endabstract

\thanks
The author was supported in part by a NUF--NAL award from the Nuffield
Foundation.
\endthanks

\endtopmatter

\centerline{\bf To appear in the Journal of Algebra}

\head Introduction \endhead
The purpose of this paper is to exhibit some remarkable links between
(i) the planar algebras on $1$-boxes arising in the work of Jones
\cite{{\bf 11}}, (ii) tabular algebras, as introduced by the author in
\cite{{\bf 7}} and (iii) the canonical bases for Hecke algebra quotients
that were defined by the author and J. Losonczy in \cite{{\bf 8}}.

In the work of Jones, a construction is given 
\cite{{\bf 11}, Example 2.2} for a kind of wreath product, $P_n^A$, of 
an associative algebra $A$
with the Temperley--Lieb algebra $TL(n, \d)$.  The construction
depends on the algebra $A$ having identity and being equipped with a
trace function.


Table algebras, which were introduced in \cite{{\bf 1}} and generalized 
in \cite{{\bf 2}}, are associative algebras with identity equipped with
distinguished bases and a natural trace function; the variant we use here 
is essentially
the same as that of Sunder's discrete hypergroups in \cite{{\bf 18}}.  
We are interested
here in the ``wreath product'' of a Temperley--Lieb algebra with a table
algebra, $A$.  The basis of $A$ equips the resulting algebra
with a ``canonical'' basis.  More precisely, Theorem 3.2.3 shows that 
this new basis is a tabular basis in the sense of \cite{{\bf 7}}.


We are particularly interested in the case where the table algebra $A$
in question is the Verlinde algebra $V_r$, which first arose in
conformal field theory \cite{{\bf 19}}.  In this case, the tabular basis of 
algebra $P_n^A$
is denoted by $\BB(n, r)$, where $n$ and $r$ are arbitrary positive
integers.  The algebra $P_n^A$ has a rich subalgebra structure
compatible with the bases $\BB(n, r)$; in particular, $\BB(n, 1)$
agrees with the usual basis of the Temperley--Lieb algebra and 
$\BB(1, r)$ agrees with the usual basis of the Verlinde algebra.  When
the parameter $n$ is allowed to vary, this produces examples of Jones' 
planar algebras, but usually we will not make this explicit and fix
$n$ in the results of this paper.

We consider pairs $(P', \BB')$ where $\BB'$ is a subset of
$\BB(n, r)$ for some $n$ and $r$ and $P'$ is a subalgebra of $P_n^A$
spanned by $\BB'$.  What is interesting (Theorem 4.2.5) is that there are
important cases of this form where $\BB'$ is precisely the
canonical basis (in the sense of \cite{{\bf 8}}) of a certain quotient of 
a Hecke algebra of type $A$, $B$, $H$ or $I$.  We conclude that there
is a close relationship between the Kazhdan--Lusztig bases of
\cite{{\bf 13}} on the one hand and certain wreath products of discrete 
hypergroups with Temperley--Lieb algebras on the other.


Although our main results consider the case where the hypergroups of
the title are Verlinde algebras, we develop the theory more 
generally because it is useful in other contexts, such as the recent
work of Rui and Xi on cyclotomic Temperley--Lieb algebras \cite{{\bf 17}}.

Only a very small part of the machinery of planar algebras is
necessary for our purposes.  In particular, we do not consider 
$C^*$-algebras and our main motivation
comes from the theory of Hecke algebras of Coxeter groups. 
We choose to use Jones' formalism because
it makes it easy to generalize Martin and Saleur's diagrammatic
setting for the blob algebra \cite{{\bf 16}} in a rigorous way.


\head 1. Table algebras \endhead

In \S1, we recall the definition of a table algebra and show that the
Verlinde algebra satisfies the table algebra axioms.  

\subhead 1.1 Definitions \endsubhead

Table algebras were defined by Arad and Blau \cite{{\bf 1}} in the
finite-dimensional, commutative case.  We tend to follow the notation of
the paper \cite{{\bf 1}}, although our definition includes some non-commutative
algebras as in \cite{{\bf 2}} and, potentially, infinite dimensional algebras.
Our algebras are more or less the same as Sunder's discrete
hypergroups \cite{{\bf 18}, \S IV}.

\definition{Definition 1.1.1}
A table algebra is a pair $(A, \BB)$, where $A$ is an associative 
unital $R$-algebra for some $\zed \leq R \leq \complex$ 
and $\BB = \{b_i : i \in I\}$ is a distinguished basis for
$A$ such that $1 \in \BB$, satisfying the following three axioms:

\item{(T1)}{The structure constants of $A$ with respect to the basis
$\BB$ lie in $\real^+$, the nonnegative real numbers.}
\item{(T2)}{There is an algebra anti-automorphism $\bar{\ }$ of $A$ whose
square is the identity and that has the property 
that $b_i \in \BB \Rightarrow \overline{b_i} \in
\BB$.  (We define $\overline{i}$ by the condition $\overline{b_i} =
b_{\bar{i}}$.)}
\item{(T3)}{Let $\k(b_i, a)$ be the coefficient of $b_i$ in $a \in A$.
Then there is a function $g: \BB \times \BB \ra \real^+$ satisfying $$
\k(b_m, b_i b_j) = g(b_i, b_m) \k(b_i, b_m \overline{b_j})
,$$ where $g(b_i, b_m)$ is independent of $j$, for all $i, j, m$.
}
\enddefinition

Following \cite{{\bf 7}}, we make the following definition which differs
slightly from the Arad--Blau notion of a ``normalized'' table algebra.

\definition{Definition 1.1.2}
A {\it normalized} table algebra $(A, \BB)$ over $R$ is one whose structure
constants lie in $\zed$ and for which the function $g$ in axiom (T3)
sends all pairs of basis elements to $1 \in \real$.  All table
algebras from now on will be normalized.
\enddefinition

\definition{Definition 1.1.3}
If $(A, \BB)$ is a table algebra and $a \in A$, we write $\supp(a)$
to denote the set of elements of $\BB$ which occur with nonzero
coefficient in $a$.
\enddefinition

Table algebras are equipped with a natural trace function; this is a
simple consequence of axiom (T3) (see \cite{{\bf 2}, \S1}).

\proclaim{Proposition 1.1.4}
Let $(A, \BB)$ be a normalized table algebra.  The linear function $t$ 
sending $a \in A$ to $\k(1, a)$ satisfies $t(xy) = t(yx)$ for all $a \in A$.
\qed\endproclaim

For our purposes, we shall need to consider certain tensor powers of
table algebras.  This construction, which is an extension of the
external direct product for groups, relies on the following simple
fact.

\proclaim{Proposition 1.1.5}
Let $(A_1, \BB_1)$ and $(A_2, \BB_2)$ be normalized table algebras
over $R \leq \complex$.  Then $(A_1
\otimes A_2, \BB_1 \otimes \BB_2)$ is a normalized table algebra,
where the multiplication on $A_1 \otimes A_2$ is given by the
Kronecker product, tensor 
products are taken over $R$ and the anti-automorphism $\bar{\
}$ of $A_1 \otimes A_2$ is defined to send $b_1 \otimes b_2$ to
$\overline{b_1} \otimes \overline{b_2}$.
\endproclaim

\demo{Proof}
The identity element of $A_1 \otimes A_2$ is $1 \otimes 1$, which is
in the basis.  Axioms (T1) and (T2) are immediate.

Consider the coefficient with which a basis element $b_1 \otimes b_2$
occurs in the product $(b_1' \otimes b_2')(b_1'' \otimes b_2'')$.
It is clear that $$
\k(b_1 \otimes b_2, (b_1' \otimes b_2')(b_1'' \otimes b_2'')) =
\k(b_1, b_1' b_1'') \k(b_2, b_2' b_2'')
.$$  Since the table algebras $(A_1, \BB_1)$ and $(A_2, \BB_2)$ are 
normalized, we have $$\eqalign{
\k(b_1, b_1' b_1'') \k(b_2, b_2' b_2'')
&= \k(b_1', b_1 \overline{b_1''}) \k(b_2', b_2 \overline{b_2''}) \cr
&= \k(b_1' \otimes b_2', (b_1 \otimes b_2)(\overline{b_1''} \otimes
\overline{b_2''})). \cr
}$$  This proves axiom (T3) and shows that the resulting table algebra
is normalized.
\qed\enddemo

\subhead 1.2 The Verlinde algebra \endsubhead

\definition{Definition 1.2.1}
Let $\{U_n(x)\}_{n \in \enn}$ be the sequence of polynomials defined 
by the conditions $U_0(x) = 1$, $U_1(x) = x$ and the recurrence relation 
$U_{n+1}(x) = x U_n(x) - U_{n-1}(x)$ for $n > 1$.  
\enddefinition

The polynomials $U_n(2x)$ are sometimes called ``type II Chebyshev
polynomials''.  We use these polynomials to define the Verlinde
algebra, which first appeared in \cite{{\bf 19}}.

\definition{Definition 1.2.2}
Let $r \geq 1$.  The Verlinde algebra, $V_r$, is defined to be the quotient of
$\zed[x]$ by the ideal generated by $U_r(x)$.  It has rank $r$, and
we equip it with a $\zed$-basis consisting of the images $u_i(x)$ of the
elements $U_i(x)$ for $0 \leq i < r$.
\enddefinition

The following result is well-known.

\proclaim{Proposition 1.2.3}
The distinguished basis for $V_r$ given by Definition 1.2.2 gives the
Verlinde algebra the structure of a table algebra, where the
automorphism $\bar{\ }$ is the identity map.
\endproclaim

\demo{Proof}
The structure constants of $V_r$  can be computed by the
Clebsch--Gordan rule.  Suppose $0 \leq n \leq n' < r$.
Then the Clebsch--Gordan rule gives, in our notation, $$
u_n(x) u_{n'}(x) = \sum_{i = 0}^{\min(n, r-n'-1)} u_{n'-n+2i}(x)
.$$  Axioms (T1) and (T2) are now satisfied.  
Note that $u_0(x)$ occurs with coefficient $1$ in this product if
$n = n'$, and with coefficient $0$ otherwise.  Axiom (T3) follows from
this observation and \cite{{\bf 1}, Lemma 2.1}.
\qed\enddemo

We conclude this section with some technical lemmas which will be useful later.

\proclaim{Lemma 1.2.4}
Denote by $w$ the element $u_{r-1}(x)$ of $V_r$.  Then, for any $0
\leq i < r$, we have $u_i(x)w = u_{r-1-i}(x)$.  In particular, $w^2 = 1$.
\endproclaim

\demo{Proof}
This is immediate from the formula in the proof of Proposition 1.2.3,
with $n' = r-1$.
\qed\enddemo

\proclaim{Lemma 1.2.5}
Consider the algebra $V_3$ over $\complex$, with basis $1 = u_0(x)$,
$y = u_1(x)$ and $z = u_2(x)$, and the algebra $V_2$ over $\complex$,
with basis $\{1', z'\}$.  Then there is a homomorphism $\phi : V_3 \ra
V_2$ such that $\phi(1) = 1'$, $\ds{\phi(y) = {{z' + 1'} \over \sqrt{2}}}$ and
$\phi(z) = z'$.
\endproclaim

\demo{Proof}
We verify that $\phi$ respects the relations $y^2 = 1+z$, $yz = zy =
y$ and $z^2 = 1$, which presents no problems.
\qed\enddemo

\head 2. Planar algebras on $1$-boxes \endhead

In \S2, we sketch Jones' construction of planar algebras on $1$-boxes arising
from associative algebras $A$, where $A$ has identity and is equipped
with a trace map.  Although this construction can be made completely rigorous,
this would take too much space so we refer the
reader to \cite{{\bf 11}} for the foundations behind the definitions.

\subhead 2.1 The Temperley--Lieb algebra \endsubhead

The concept of a $k$-box, for $k \in \enn$, is defined in \cite{{\bf 11},
Definition 1.1}, as follows.

\definition{Definition 2.1.1}
Let $k$ be a nonnegative integer.  The {\it standard $k$-box}, ${\Cal
B}_k$, is the set $\{(x, y) \in \real^2 : 0 \leq x \leq k + 1, \ 0
\leq y \leq 1\}$, together with the $2k$ marked points $$\eqalign{
&1 = (1, 1), \ 
2 = (2, 1), \ 
3 = (3, 1), \ 
\ldots, \ k = (k, 1), \cr
&k + 1 = (k, 0), \
k + 2 = (k-1, 0), \
\ldots, \ 
2k = (1, 0).\cr
}$$
\enddefinition

We summarise the definition of the algebra $\Pl_k(\emptyset)$ from
\cite{{\bf 11}}.

\definition{Definition 2.1.2}
Let $k$ be a nonnegative integer.  An element of $T_k(\emptyset)$
consists of a finite number of oriented disjoint curves (which we
usually call ``edges''), smoothly
embedded in the standard $k$-box, under smooth orientation-preserving
diffeomorphisms of $\real^2$.  Curves may be closed (isotopic to
circles) but not if their endpoints coincide with marked points of the
box.  The marked points of the box are endpoints of curves, which meet the
box transversely.  Otherwise, the curves are disjoint from the box.  
The orientations of the curves must satisfy the following two conditions.

\item{(i)}{A curve meeting the $r$-th marked point of the standard
$k$-box, where $r$ is odd, must exit the box at that point.}

\item{(ii)}{Each connected component of the complement of the union of
the curves in the standard $k$-box may be oriented in such a way that
the orientation of a curve coincides with the orientation induced as
part of the boundary of the connected component.}
\enddefinition

\example{Example 2.1.3}
Let $k = 8$.  An element of $T_8(\emptyset)$ is shown in Figure 1.
Note that there are 10 connected components as in Definition 2.1.2
(ii), of which precisely 7 inherit a clockwise orientation.
\endexample

\topcaption{Figure 1} Typical element of $T_8(\emptyset)$ \endcaption
\centerline{
\hbox to 4.527in{
\vbox to 1.375in{\vfill
        \includegraphics{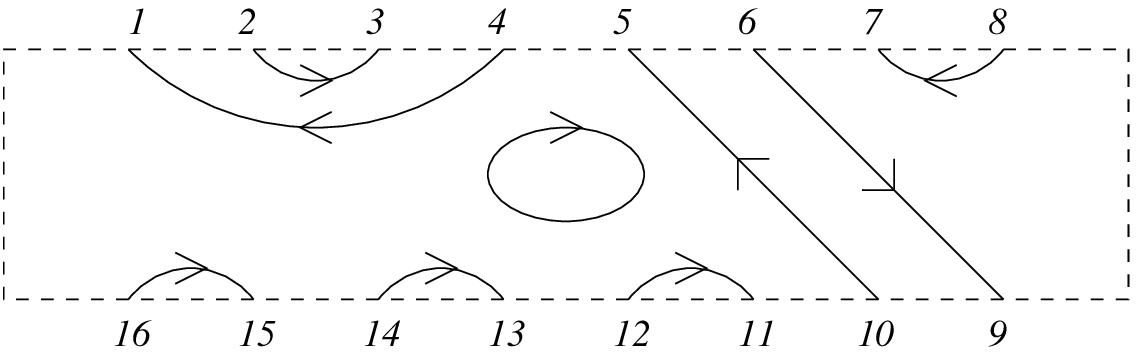}
}
\hfill}
}

The following definition is a special case of \cite{{\bf 11}, Definition
1.8}.

\definition{Definition 2.1.4}
Let $k$ be a nonnegative integer.
The associative algebra $\Pl_k(\emptyset)$ over a commutative ring $R$
with $1$ is the free $R$-module having $T_k(\emptyset)$ as a basis,
with multiplication defined as follows.  If $T_1, T_2 \in
T_k(\emptyset)$, the product $T_1T_2$ is the element of
$T_k(\emptyset)$ obtained by placing $T_1$ on top of $T_2$,
rescaling vertically by a scalar factor of $1/2$ and applying the
appropriate translation to recover a standard $k$-box.

Note that we may assume that the resulting curves are smooth.  The
orientations will match up automatically.
\enddefinition

The Temperley--Lieb algebra can be easily defined in terms of this
formalism, as shown in \cite{{\bf 11}, Definition 2.1}.

\definition{Definition 2.1.5}
Let $R$ be a commutative ring with $1$.  The Temperley--Lieb algebra,
$TL(n, \d)$, is the free $R[\d]$-module with basis given by the elements
of $T_n(\emptyset)$ with no closed loops.  The multiplication is
inherited from the multiplication on $\Pl_n(\emptyset)$ except that
one multiplies by a factor of $\d$ for each resulting closed loop and
then discards the loop.
\enddefinition

We usually consider $TL(n, \d)$ to be an algebra defined over $\A :=
\zed[v, v^{-1}]$, where $\d = v + v^{-1}$.  The Laurent polynomial $v
+ v^{-1}$ is often denoted by $[2]$, which will be our preferred notation.

\subhead 2.2 Planar algebras on $1$-boxes \endsubhead

In \S2.2, we recall from \cite{{\bf 11}, Example 2.2} the construction of
the algebra $P_n^A$ from the Temperley--Lieb algebra $TL(n, \d)$ and
the associative $R$-algebra $A$, where $R$ is a commutative ring
containing $\d$.  
The algebra $A$ is assumed to have identity and a trace functional 
$\tr: A \ra R$ with $\tr(ab) = \tr(ba)$ and $\tr(1) = \d$.

\definition{Definition 2.2.1}
Let $A$ be as above, and let $k$ be a nonnegative integer.
We define the tangles $T_k(A)$ to be those that arise from elements
of $T_k(\emptyset)$ by adding zero or more $1$-boxes labelled by
elements of $A$ to each edge.  An edge of $T_k(A)$ that is not a loop
is called {\it propagating} if its endpoints have different
$y$-values, and {\it non-propagating} otherwise.  An edge of $T_k(A)$
which is not a loop is called
{\it transitional} if its endpoints lie on different sides of the line
$x = 3/2$, and {\it non-transitional} otherwise.  Transitional edges
may or may not be propagating.
\enddefinition

Figure 2 shows a typical element of $T_8(A)$
which $a, b, c, d, e$ are some elements of the algebra $A$.  There are
2 propagating edges and 6 non-propagating edges.  There are 2
transitional edges (those emerging from points $4$ and $16$) and 6
non-transitional edges.

\topcaption{Figure 2} Typical element of $T_8(A)$ \endcaption
\centerline{
\hbox to 4.527in{
\vbox to 1.375in{\vfill
        \includegraphics{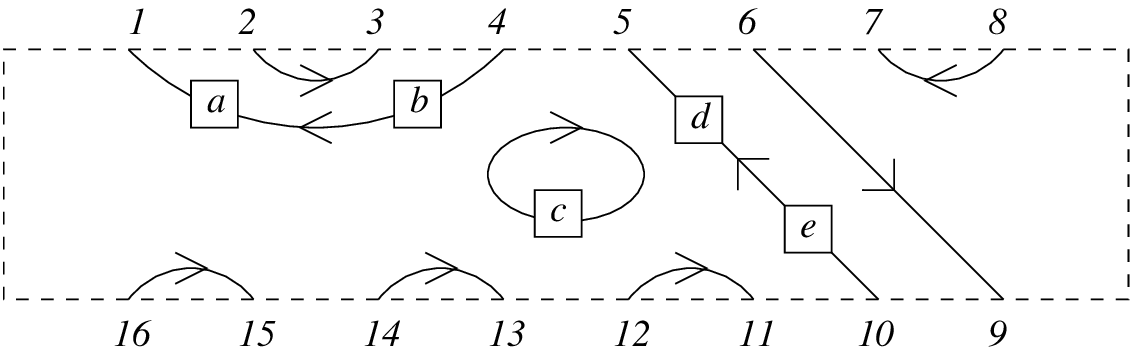}
}
\hfill}
}

Definitions 2.1.4 and 2.1.5 generalize naturally to this situation, as
follows.

\definition{Definition 2.2.2}
Let $k$ be a nonnegative integer and let $A$ be an $R$-algebra (as
before) with a free $R$-basis, $\{a_i : i \in I\}$, where $1 \in \{a_i\}$.
The associative $R$-algebra $P_k^A$
is the free $R$-module having as a basis those elements of $T_k(A)$ 
satisfying the conditions that
\item{(i)}{all labels on edges are basis elements $a_i$,}
\item{(ii)}{each edge has precisely one label and}
\item{(iii)}{there are no closed loops}.

The multiplication is defined in the case where $T_1$ and $T_2$ are
basis elements of $P_k^A$ as above, and extended bilinearly.
To calculate the product $T_1T_2$, place $T_1$ on top of $T_2$,
rescale vertically by a scalar factor of $1/2$ and apply the
appropriate translation to recover a standard $k$-box.  Next, apply
relations (a), (b) and (c) below to express the product as an
$R$-linear combination of basis elements, and finally, apply relation (d)
below to remove any loops, multiplying by the scalar shown for each
loop removed.
\enddefinition

\topcaption{Figure 3} Relation (a) of Definition 2.2.2 \endcaption
\centerline{
\hbox to 2.513in{
\vbox to 1.027in{\vfill
        \includegraphics{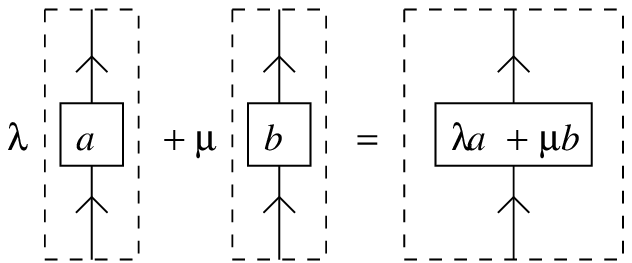}
}
\hfill}
}

\topcaption{Figure 4} Relation (b) of Definition 2.2.2 \endcaption
\centerline{
\hbox to 1.166in{
\vbox to 1.277in{\vfill
        \includegraphics{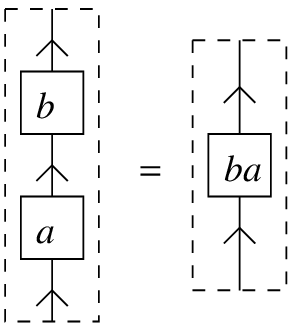}
}
\hfill}
}

\topcaption{Figure 5} Relation (c) of Definition 2.2.2 \endcaption
\centerline{
\hbox to 1.166in{
\vbox to 1.027in{\vfill
        \includegraphics{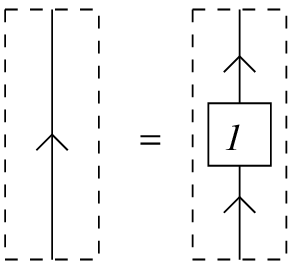}
}
\hfill}
}


\topcaption{Figure 6} Relation (d) of Definition 2.2.2 \endcaption
\centerline{
\hbox to 1.291in{
\vbox to 1.027in{\vfill
        \includegraphics{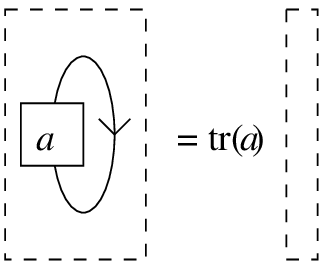}
}
\hfill}
}

\remark{Remark 2.2.3}
The direction on the arrow in relation (d) is immaterial.  For a proof
that this procedure does define an associative algebra, the reader is
referred to \cite{{\bf 11}, Example 2.2}.
\endremark

\subhead 2.3 Planar algebras from table algebras \endsubhead

If $A$ is a normalized table algebra, we can be more specific about
the construction in \S2.2.

\definition{Definition 2.3.1}
Let $A$ be a normalized table algebra over $R$ (containing $\d$)
with distinguished
basis $\BB$.  Then, for $n \in \enn$, we define the $R$-algebra $P_n^A$ to
be that arising from Definition 2.2.2 with the basis $\BB$ as the
distinguished basis; for the trace map, we take $\d . t$ where $t$ is
the trace in Proposition 1.1.4.
In this case, we refer to the $R$-basis of $P_n^A$ that arises in
Definition 2.2.2 as the {\it canonical basis}, $\BB_n^A$, of $P_n^A$.
\enddefinition

Note that if $A$ is a one-dimensional unital algebra, this
construction yields the Temperley--Lieb algebra (Definition 2.1.5).

For later purposes, we require the following natural anti-automorphism
of the algebra $P_n^A$.

\proclaim{Lemma 2.3.2}
There is an linear anti-automorphism, $*$, of $P_n^A$ permuting the canonical
basis.  The image, $b^*$, of a basis element $b$ under this map is
obtained by reflecting $b$ in the line $y = 1/2$, reversing the
direction of all the arrows and replacing each $1$-box labelled by
$b_i \in \BB$ by a $1$-box labelled by $b_{\bar{i}}$.
\endproclaim

\demo{Proof}
This follows easily from Definition 2.2.2, because the
anti-automorphism of the table algebra is $\zed$-linear and the trace of
Proposition 1.1.4 satisfies $t(x) = t(x^*)$ for all $x$.
\qed\enddemo

\proclaim{Lemma 2.3.3}
Let $A$ be a normalized table algebra as in Definition 2.3.1.
Consider a closed loop associated with an element of $P_n^A$ that
carries precisely two $1$-boxes, each of which is labelled by a basis
element.  If the labels are $b$ and $b'$, then removal of the loop
results in multiplication by zero unless $\overline{b} = b'$, in which
case removal of the loop results in multiplication by $\d$.
\endproclaim

\demo{Proof}
Let $t$ be the trace map of Proposition 1.1.4.  Then $t(b b')$ is
$1$ if $\overline{b} = b'$ and $0$ otherwise.  The result follows from
Definition 2.3.1.
\qed\enddemo

As mentioned in \cite{{\bf 11}, Example 2.2}, the algebra $P_n^A$ is closely
related to $A^{\otimes n}$.  The following result is a special case of
this relationship, and is useful for later purposes.

\proclaim{Proposition 2.3.4}
Let $A$ be a normalized table algebra over $R$ with distinguished
basis $\BB$.  There is an isomorphism $\rho$ of $R$-algebras
from $(A^{\otimes n}, \BB^{\otimes n})$ to the subalgebra of $P_n^A$ spanned by
all canonical basis elements with no non-propagating edges.  The
isomorphism takes basis elements to basis elements.
\endproclaim

\demo{Proof}
Let $b = b_{i_1} \otimes b_{i_2} \otimes \cdots \otimes b_{i_n}$ be a
typical basis element from the set $\BB^{\otimes n}$.  This element is
sent by the isomorphism, $\rho$, to a canonical basis element of $P_n^A$ with
no propagating edges, where the decoration on the $k$-th propagating
edge (counting from $1$ to $n$, starting at the left) is $b_{i_k}$ if
$k$ is odd, and $b_{\overline{i_k}}$ if $k$ is even.  (Note that the
$k$-th propagating edge points upwards if and only if $k$ is odd; this
is a consequence of the orientation on
the standard $n$-box.)  The relations (a)--(d) of Definition 2.2.2
show that $\rho$ is $R$-linear and bijective.

We now invoke an idea from Lemma 2.3.2: arrows may be reversed at the expense
of applying the table algebra automorphism to their labels.
It is now clear that $\rho$ defines an isomorphism of $R$ algebras.
\qed\enddemo

\subhead 2.4 The algebra $\D_n^A$ \endsubhead

Let $A$ be a normalized table algebra over $R$ and let $P_n^A$ be the
algebra arising from $A$ as in Definition 2.3.1.  We now define a certain
subalgebra of $P_n^A$ which turns out to be useful for our purposes.

\definition{Definition 2.4.1}
Let $A$ be a normalized table algebra and let $D \in T_n(A)$ for some
$n$ (see Figure 2).  Of the connected components of the complement of
the union of the curves of $D$ in the standard $n$-box, there is a
unique component which is bounded on the left by the line $x = 0$.  We
call this the {\it principal} connected component of $D$.
\enddefinition

\example{Example 2.4.2}
In Figure 2, there are precisely 6 edges adjacent to the principal
connected component of the tangle shown, including the closed loop.
\endexample

\definition{Definition 2.4.3}
Maintain the above notation.  A canonical basis element $b$ of $P_n^A$
is said to be {\it exposed} if the only edges labelled by nonidentity
basis elements of $A$ are adjacent to the principal connected
component of $b$.  We denote the free $R$-submodule of $P_n^A$ 
spanned by the exposed basis elements by $\D_n^A$.
\enddefinition

\proclaim{Proposition 2.4.4}
Let $A$ be a normalized table algebra over $R$ and let $P_n^A$ and
$\D_n^A$ be as above.  Then $\D_n^A$ is a subalgebra of $P_n^A$.
\endproclaim

\demo{Proof}
Let $T_1$ and $T_2$ be exposed basis elements of $P_n^A$.  It suffices
to show that the product $T_1 T_2$ is a linear combination of exposed basis
elements.  Consider an edge, $e$, in the product $T_1T_2$ which is labelled
by anything other than a sequence of $1$-boxes containing $1 \in A$.
Then $e$ must be an extension of an edge of $T_i$ (where $i = 1$ or
$2$) that is adjacent to the principal connected component of $T_i$.
It follows that $e$ is adjacent to the principal connected component
of $T_1T_2$, and thus that $e$ is exposed.
\qed\enddemo

\definition{Definition 2.4.5}
Let $n, r$ be positive integers.  Let $A$ be the Verlinde algebra, $V_r$, with
its distinguished basis as in Proposition 1.2.3.  We define $\D(n, r)
:= \D_n^A$ and $P(n, r) := P_n^A$.
\enddefinition

The algebra $P(n, r)$ has a useful involution that permutes the basis
elements, defined as follows.

\proclaim{Lemma 2.4.6}
Let $\omega : P(n, r) \ra P(n, r)$ be the linear map defined by sending
a canonical basis element $D$ to the element of $P(n, r)$
which differs from $D$ only in that any transitional edges labelled by 
$b \in V_r$ are relabelled by $u_{r-1}b$.  Then $\omega$ is an 
$R$-algebra automorphism that permutes the canonical basis and
preserves $\D(n, r)$ setwise.
\endproclaim

\demo{Proof}
By Lemma 1.2.4, $u_{r-1}b$ is indeed a basis element of $V_r$, so $\omega$
permutes the canonical basis elements.  It is clear that $\omega$ is
linear and bijective.
Also by Lemma 1.2.4, the element $w := u_{r-1}$ satisfies $w^2 = 1$.


Consider an edge or loop, $e$, in the product $DD'$,
and suppose that $e$ is labelled by $g \in V_r$ (where $g$ is not
necessarily a basis element).  Let $e'$ be the edge
or loop corresponding to $e$ in the product $\omega(DD')$, and let its label
be $g' \in V_r$.  Now $e$ may include 
either (a) an even or (b) an odd number of transitional edges from
$\omega(D)$ and $\omega(D')$.
In case (a), $e'$ crosses the line $x = 3/2$ an
even number of times and must therefore be a non-transitional
edge or a loop; furthermore, $g = g'$ as $w^2 = 1$.
In case (b), $e$ crosses the line $x = 3/2$ an odd number of
times and must therefore be
a transitional edge (and not a loop), and we have $g' = wg$.  It
follows from these observations that $\omega(DD') =
\omega(D)\omega(D')$, and thus that $\omega$ is a homomorphism.

To prove the last assertion, note that if $D$ has any 
(\idest $2$) transitional edges, they must be exposed, as they cross
the line $x = 3/2$.  This implies that the subalgebra $\D(n, r)$ is
invariant under $\omega$.
\qed\enddemo

\remark{Remark 2.4.7}
An interesting combinatorial problem arising from these definitions 
might be to determine an explicit formula for the rank of the algebra 
$\D(n, r)$ for fixed values of $r$.  For $r = 1$, the rank
is the $n$-th Catalan number, $\ds{{1 \over {n+1}} {{2n} \choose
n}}$; for $r = 2$, the rank is $\ds{{2n} \choose n}$.
For $r = 3$, the sequence appears to be related to the number of
chains in rooted plane trees (sequence $w_3$ in \cite{{\bf 14}, \S2}).
\endremark

\head 3. Tabular algebras \endhead


We recall from \cite{{\bf 7}, \S1.3} the definition of a tabular algebra.
The example to keep in mind throughout
is the algebra $P_n^A$, which will be seen in Theorem 3.2.3 to be a
tabular algebra with basis $\BB_n^A$ satisfying the five axioms below.

\subhead 3.1 Tabular algebras and sub-tabular algebras \endsubhead

\definition{Definition 3.1.1}
Let $\A = \zed[v, v^{-1}]$.  A {\it tabular algebra} is an
$\A$-algebra $A$, together with a table datum 
$(\Lambda, \Gamma, B, M, C, *)$ where:

\item{(A1)}
{$\Lambda$ is a finite poset.  For each $\l \in \Lambda$, 
$(\Gamma(\l), B(\l))$ is
a normalized table algebra over $\zed$ and $M(\l)$ is a finite set.
The map $$
C : \coprod_{\l \in \Lambda} \left( M(\l) \times B(\l) \times M(\l)
\right) \rightarrow A
$$ is injective with image an $\A$-basis of $A$.  We assume
that $\Im(C)$ contains a set of mutually orthogonal idempotents 
$\{1_\e : \e \in {\Cal E}\}$ such that 
$A = \sum_{\e, \e' \in {\Cal E}} (1_\e A 1_{\e'})$ and such that for each
$X \in \Im(C)$, we have $X = 1_\e X 1_{\e'}$ for some $\e, \e' \in
{\Cal E}$.
A basis arising in
this way is called a {\it tabular basis}.  
}
\item{(A2)}
{If $\l \in \Lambda$, $S, T \in M(\l)$ and $b \in B(\l)$, we write
$C(S, b, T) = C_{S, T}^{b} \in A$.  
Then $*$ is an $\A$-linear involutory anti-automorphism 
of $A$ such that
$(C_{S, T}^{b})^* = C_{T, S}^{\overline{b}}$, where $\bar{\ }$ is the
table algebra anti-automorphism of $(\Gamma(\l), B(\l))$.
If $g \in \complex(v) \otimes_\zed \Gamma(\l)$ is such that 
$g = \sum_{b_i \in B(\l)} c_i b_i$ for some scalars $c_i$ 
(possibly involving $v$), we write
$C_{S, T}^g \in \complex(v)\otimes_\A A$ 
as shorthand for $\sum_{b_i \in B(\l)} c_i C_{S, T}^{b_i}$.  We write
$\bc_\l$ for the image under $C$ of $M(\l) \times B(\l) \times M(\l)$.}
\item{(A3)}
{If $\l \in \Lambda$, $g \in \Gamma(\l)$ and $S, T \in M(\l)$ then for all 
$a \in A$ we have $$
a . C_{S, T}^{g} \equiv \sum_{S' \in M(\l)} C_{S', T}^{r_a(S', S) g}
\mod A(< \l),
$$ where  $r_a (S', S) \in \Gamma(\l)[v, v^{-1}] = \A \otimes_\zed
\Gamma(\l)$ is independent of $T$ and of $g$ and $A(< \l)$ is the
$\A$-submodule of $A$ generated by the set $\bigcup_{\mu < \l} \bc_\mu$.}
\enddefinition

Next, we recall the $\ba$-function associated to a tabular algebra $A$.

\definition{Definition 3.1.2}
Let $g_{X, Y, Z} \in \A$ be one of the structure constants for the tabular
basis $\Im(C)$ of $A$, namely $$
X Y = \sum_Z g_{X, Y, Z} Z
,$$ where $X, Y, Z \in \Im(C)$.   Define, for $Z \in \Im(C)$, $$
\afn(Z) = 
\max_{X, Y \in \Im(C)} \deg(g_{X, Y, Z})
,$$ where the degree of a Laurent polynomial is taken to be 
the highest power of $v$ occurring with nonzero coefficient.  We
define $\g_{X, Y, Z} \in \zed$ to be the coefficient of $v^{\afn(Z)}$ in $g_{X,
Y, Z}$; this will be zero if the bound is not achieved.
\enddefinition

Using the notion of $\ba$-function, we recall the definition of ``tabular
algebras with trace''.

\definition{Definition 3.1.3}
A {\it tabular algebra with trace} is a tabular algebra in the sense
of Definition 3.1.1 that satisfies the conditions (A4) and (A5) below.

\item{(A4)}{Let $K = C_{S, T}^b$, $K' = C_{U, V}^{b'}$ and 
$K'' = C_{X, Y}^{b''}$ lie in $\Im(C)$.  Then the
maximum bound for $\deg(g_{K, K', K''})$
in Definition 3.1.2 is achieved if and only if $X = S$, $T = U$, $Y =
V$ and $b'' \in
\supp(bb')$ (see Definition 1.1.3).  If these conditions all hold and
furthermore $b = b' = b'' = 1$, we require $\g_{K, K', K''} = 1$.}
\item{(A5)}{There exists an $\A$-linear function $\t : A \ra \A$
(the {\it tabular trace}), such that $\t(x) = \t(x^*)$ for all $x \in
A$ and $\t(xy) = \t(yx)$ for all $x, y \in A$, that has the 
property that for every
$\l \in \Lambda$, $S, T \in M(\l)$, $b \in B(\l)$ and $X = C_{S,
T}^b$, we have $$
\t(v^{\afn(X)} X) = 
\cases 1 \mod v^{-1} \A^- & \text{ if } S = T \text{ and } b = 1,\cr
0 \mod v^{-1} \A^- & \text{ otherwise.} \cr
\endcases
$$  Here, $\A^- := \zed[v^{-1}]$.}
\enddefinition

The main results of this paper may be described in terms of ``sub-tabular
algebras'', which we now introduce.

\definition{Definition 3.1.4}
Let $(A, \BB)$ be a tabular algebra $A$ together with its tabular
basis $\BB = \Im(C)$.  Then a sub-tabular algebra is a pair $(A',
\BB')$, where $\BB' \subseteq \BB$ and $A'$ is a subalgebra of $A$
with $\BB'$ as a basis.
\enddefinition

Note that a sub-tabular algebra may or may not be a tabular algebra in
its own right.  Some, but not all, of the examples of sub-tabular
algebras that we consider in \S4.2 are also tabular algebras.

\subhead 3.2 Tabular structure of $P_n^A$ \endsubhead

In \S3.2, we suppose $(A, \BB)$ is a table algebra and we show that
$P_n^A$, equipped with its canonical basis $\BB_n^A$, is a tabular 
algebra in a natural way.  This involves equipping the algebra $P_n^A$
with a certain trace, which is easily done using the formalism of
planar algebras.

Recall from \cite{{\bf 11}, Definition 1.2.8} that for a spherical planar
algebra (\idest one for which relations such as that in Figure 6 may
be transformed into their mirror images) we may define a trace on a
$k$-box, $x$, by the following procedure.  First, for each $1 \leq
i \leq k$, join the point $i$ to the point $2k + 1 - i$ using 
suitably oriented non-intersecting curves, and then repeatedly apply
relation (d) of Definition 2.2.2 (with respect to a suitable trace on
the underlying table algebra) to obtain $\tr(x)$.

\definition{Definition 3.2.1}
Let $(A, \BB)$ be a table algebra over $\A = \zed[v, v^{-1}]$ and let 
$\d = v + v^{-1}$.
We define the trace $\t : P_n^A \ra R$ by $\t(x) := v^{-n} \tr(x)$,
where $\tr = \tr_L = \tr_R$ is Jones' trace from \cite{{\bf 11},
Definition 1.28} compatible with the trace $\d . t$ on $A$, with
$t$ as in Proposition 1.1.4.
\enddefinition

\example{Example 3.2.2}
The identity element $1 \in P_n^A$ satisfies $\tr(1) = \d^n$ and
$\t(1) = (1 + v^{-2})^n$.  

The element $x$ of $P_8(A)$ shown in Figure
2 satisfies $\tr(x) = \d^3 t(ab) t(c) t(de)$ and $\t(x) = (v^{-5} + 3
v^{-7} + 3 v^{-9} + v^{-11}) t(ab) t(c) t(de)$.
\endexample

Although its proof is similar to \cite{{\bf 7}, Theorem 5.2.5},
the next result is much more general.  Its relevance will become
clearer in \S4.3.

\proclaim{Theorem 3.2.3}
The algebra $P_n^A$ equipped with its canonical basis $\BB_n^A$ and 
the trace $\t$ of Definition 3.2.1 is a tabular algebra with trace.
\endproclaim

\demo{Proof}
Let $\Lambda$ be the set of integers $r$ with $0
\leq r \leq n$ and $n - r$ even, ordered in the usual way.  

For $\l \in \Lambda$, let $(\G(\l), B(\l))$ be the $\l$-th tensor
power of the table algebra $(A, \BB)$ with the basis and
anti-automorphism induced by Proposition 1.1.5.


Let $M(\l)$ be the set of possible configurations of non-propagating
edges with endpoints on the line $y = 1$ that arise from an element of 
$\BB_n^A$.  (Note that, in this case, the number of non-propagating
edges involved will be $(n - \l)/2$.)
Let $b = b_{i_1} \otimes b_{i_2} \otimes \cdots \otimes
b_{i_\l}$ be a
basis element of $B(\l)$ and let $m$ and $m'$ be elements of $M(\l)$.
The map $C$ produces a basis element in $\BB_n^A$ from the triple 
$(m, b, m')$ as follows.  Turn the half-diagram corresponding to $m'$
upside down, reverse the directions of all the arrows and relabel all
$1$-boxes labelled by $b_i \in \BB$ so they are labelled by
$b_{\bar{i}}$.
Join any free marked points in the line $y = 0$ free marked points in 
the line $y = 1$ so that they do not intersect.  Orient any new edges
according to the orientation of the standard $n$-box.  Decorate the $\l$
propagating edges with the basis element $b$ exactly as in the proof
of Proposition 2.3.4.  (See Example 3.2.4 below for an illustration.)

The map $*$ is as given in Lemma 2.3.2.

It is clear that the image of $C$ as above is the canonical basis of
$P_n^A$.  Furthermore, $C$ contains the identity basis element, so
axiom (A1) holds.  Axiom (A2) follows from Lemma 2.3.2.  Axiom (A3)
follows from consideration of a product of a basis element $a$ with a
basis element $C_{S, T}^b \in \bc_\l$ in the case where $a C_{S, T}^b$ does
not lie in $P_n^A(< \l)$.  In this case, the structure constants occurring
are essentially unaffected by changing the
configuration of non-propagating edges at the bottom end of the basis element
$C_{S, T}^b$.  The part of axiom (A3) guaranteeing independence from
$g \in \Gamma(\l)$ follows from the associativity of the 
algebra structure on $(A^{\otimes
n}, \BB^{\otimes n})$ given by Proposition 2.3.4.

For axiom (A4), we claim that if $D \in \bc_\l$, we have 
$\afn(D) = a'(D) := (n - \l)/2$, \idest the
$\afn$-function evaluated at a diagram is half the number of
non-propagating edges in that diagram.  Let $D = C_{S, T}^b \in
\bc_\l$.  By Lemma 2.3.3, $C_{S, S}^1 D = [2]^{a'(D)} 
C_{S, T}^b$, so $\afn(D) \geq a'(D)$.  Conversely, the diagram
calculus shows that if $D'$ and $D''$ are canonical basis elements for
$P_n^A$, the number of loops formed in the product $D'D''$ is
bounded above both by $a'(D')$ and $a'(D'')$; this implies that the
structure constants appearing in $D'D''$ have degree bounded in the
same way.  Since $D$ can only appear
in a product $D'D''$ if $a'(D') \leq a'(D)$ and $a'(D'') \leq a'(D)$, we
have $\afn(D) \leq a'(D)$.  The claim follows.

The above argument also implies that the only way the
$\afn$-function bound can be achieved is if the three basis elements
$D', D'', D$ concerned come from the same $\bc_\l$.  
The statement of Lemma 2.3.3 shows that the
bound can only be achieved if the pattern of edges at the bottom of
$D'$ is the same as the pattern of edges at the top of $D''$, except
that all the directions on the edges have been reversed and all the
labels $b_i$ have been changed to $b_{\bar{i}}$.  In this
case, we may set $D' = C_{S, T}^b$ and $D'' = C_{T, U}^{b'}$, and
properties of the diagram calculus give $D'D'' = [2]^{\afn(D)} C_{S,
U}^{bb'}$.  The assertions of axiom (A4) all follow easily.

Finally, we prove axiom (A5).  Consider a basis element $D$.  It is
clear by symmetry of the definitions that $\t(D) = \t(D^*)$, and thus
that $\t(x) = \t(x^*)$ for all $x \in P_n^A$.  To prove the other
requirements of the axiom,
we note that the diagram corresponding to $D$
has $2k$ non-propagating edges and $r$ propagating
edges, where $2k + r = n$ and $k = \afn(D)$. 
To calculate $\t(D)$, we join each point $i$ to point $2n + 1 - i$ as
described before Definition 3.2.1.  An elementary analysis shows that
the number of loops formed is at most $k + r = n - \afn(D)$, and also that
this bound is achieved only 
if whenever point $i$ is connected to $j$, we must have point $2n
+ 1 - i$ connected to point $2n + 1 - j$.  If there are indeed $k +
r$ loops and $\t(D) \ne 0$, we also require $D = D^*$ (or one of the 
loops would give trace zero by Lemma 2.3.3) and we require all
non-propagating edges to be labelled by $1$ (or one of the loops would
contain exactly one non-identity box and give trace zero); in other
words, $D = C_{S, S}^1$ for some $S$.  If the
above bound on loops is achieved and $D = C_{S, S}^1$ for some $S$, we have 
$\t(D) = v^{-n} (v + v^{-1})^{n - \afn(D)}$ and it
follows that $\t(v^{\afn(D)}D) = 1 \mod v^{-1} \A^-$ as required.  In
the other cases, the bound is not achieved or $\t(D) = 0$, and we
have $\t(v^{\afn(D)}D) = 0 \mod v^{-1} \A^-$.
Axiom (A5) follows, completing the proof.
\qed\enddemo

\example{Example 3.2.4}
Suppose $A$ is a normalized table algebra of rank 3, with $\BB = \{1,
g, h\}$ and $\bar{g} = h$.  (For example, we could take $\BB = \zed_3$
and $\bar{\ }$ to be inversion.)
Let $n = 8$ and $\l = 2$.  Let $m, m' \in M(2)$ be as shown in Figures
7 and 8 respectively, and let $b = h \otimes h \in \BB^{\otimes 2}$.
Then the element $C_{m, m'}^b$ is as shown in Figure 9.
\endexample

\topcaption{Figure 7} The element $m \in M(2)$ of Example 3.2.4 \endcaption
\centerline{
\hbox to 3.541in{
\vbox to 0.541in{\vfill
        \includegraphics{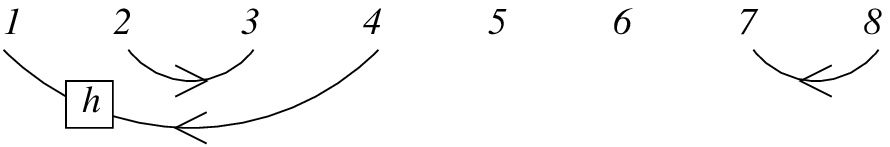}
}
\hfill}
}

\topcaption{Figure 8} The element $m' \in M(2)$ of Example 3.2.4 \endcaption
\centerline{
\hbox to 3.541in{
\vbox to 0.541in{\vfill
        \includegraphics{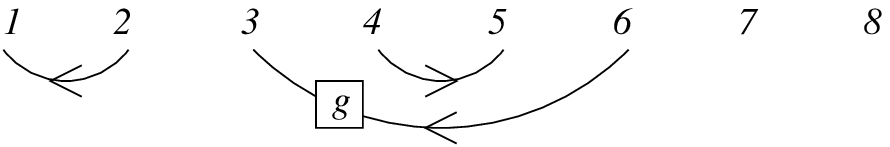}
}
\hfill}
}
\topcaption{Figure 9} The element $C(m, b, m')$ of Example 3.2.4 \endcaption
\centerline{
\hbox to 4.527in{
\vbox to 1.375in{\vfill
        \includegraphics{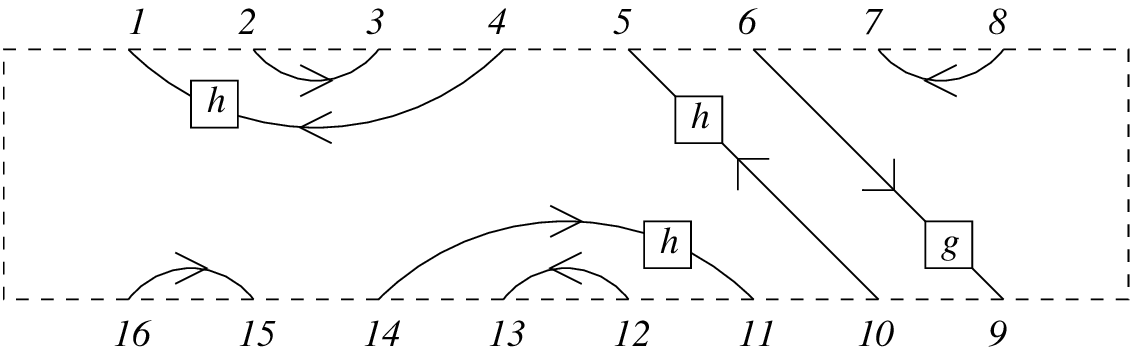}
}
\hfill}
}

\remark{Remark 3.2.5}
The algebra $\D_n^A$ is also a tabular algebra with trace.  A suitable
table datum may be obtained from the table datum for $P_n^A$ by
taking $\Lambda$ as for $P_n^A$,
considering suitable subsets of $M(\l)$, $\Gamma(\l)$ and $B(\l)$ for
each $\l \in \Lambda$, and
restricting $C$, $*$ and $\t$ to the appropriate domains.  We
leave the details to the reader.
\endremark

\head 4. Relationship with Hecke algebra quotients \endhead

In \S4.1, we recall the definition of certain Hecke algebra quotients
which we call generalized Temperley--Lieb algebras, and the
construction of their canonical bases.  These bases are related, via the
theory of tabular algebras, to the algebras $\D(n, r)$, as we explain
in the main results of \S4.2.

\subhead 4.1 Generalized Temperley--Lieb algebras and canonical bases 
\endsubhead

Let $X$ be a Coxeter graph, of arbitrary type,
and let $W(X)$ be the associated Coxeter group with distinguished
set of generating involutions $S(X)$.  Denote by $\H(X)$ the Hecke
algebra associated to $W(X)$.  Let $\A = \zed[v, v^{-1}]$ as usual.  
The $\A$-algebra $\H(X)$ has 
a basis consisting of elements $T_w$, with $w$ ranging over $W(X)$, 
that satisfy $$T_s T_w = 
\cases
T_{sw} & \text{ if } \ell(sw) > \ell(w),\cr
q T_{sw} + (q-1) T_w & \text{ if } \ell(sw) < \ell(w),\cr
\endcases$$ where $\ell$ is the length function on the Coxeter group
$W(X)$, $w \in W(X)$, and $s \in S(X)$.
The parameter $q$ is equal to $v^2$.


Let $J(X)$ be the two-sided ideal of $\H(X)$ generated by the elements $$
\sum_{w \in \lan s, s' \ran} T_w,
$$ where $(s, s')$ runs over all pairs of elements of $S(X)$
that correspond to adjacent nodes in the Coxeter graph, and $\lan s,
s' \ran$ is the group generated by the pair $(s, s')$.
(If the nodes corresponding to $(s, s')$ are connected by a
bond of infinite strength, then we omit the corresponding relation.)

\definition{Definition 4.1.1}
Following Graham \cite{{\bf 4}, Definition 6.1}, we define the generalized
Temperley--Lieb algebra $TL(X)$ to be
the quotient $\A$-algebra $\H(X)/J(X)$.  We denote the corresponding
epimorphism of algebras by $\th : \H(X) \ra TL(X)$.
\enddefinition

\topinsert
\topcaption{Figure 10} Coxeter graphs $X$ corresponding to finite rank
algebras $TL(X)$ \endcaption
\vskip 20pt
\centerline{
\hbox to 3.319in{
\vbox to 4.222in{\vfill
        \includegraphics{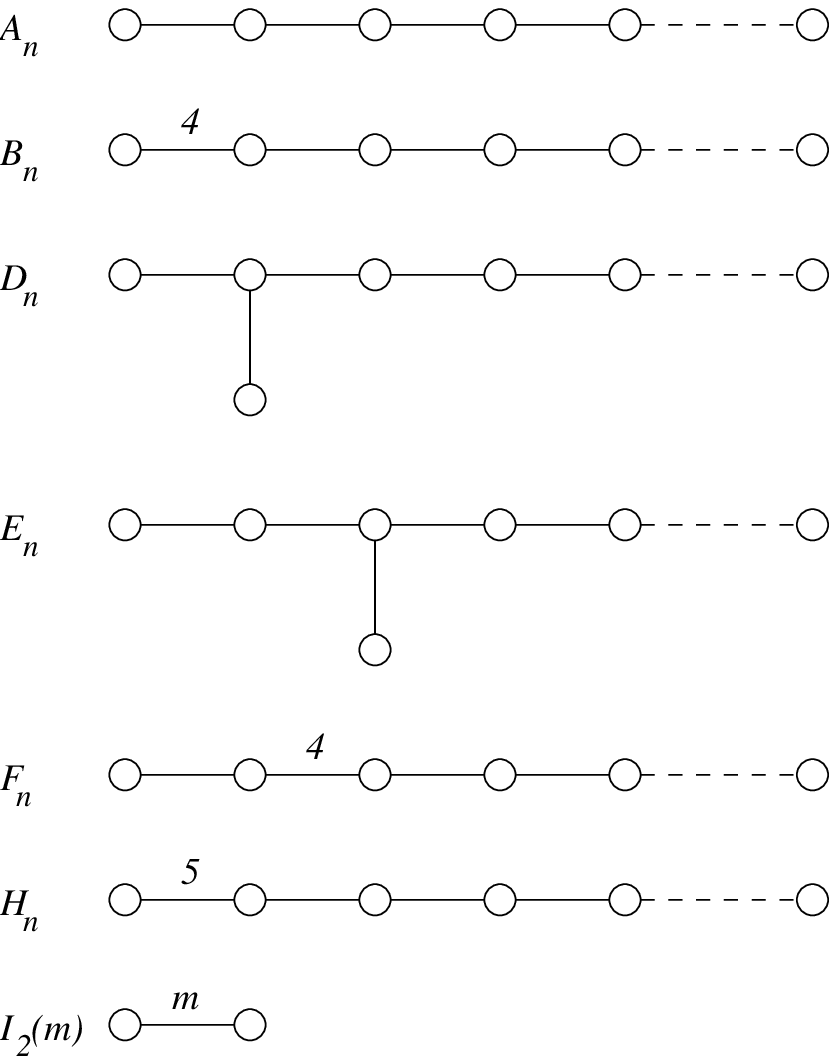}
}
\hfill}
}
\endinsert

The algebra $TL(X)$ may be of finite or infinite rank, and may be
of finite rank even when it is the quotient of a Hecke algebra of
infinite rank, as the following result shows.

\proclaim{Proposition 4.1.2 (Graham)}
The algebra $TL(X)$ corresponding to a Coxeter graph $X$ is of finite
rank if and only if its graph is one of those appearing in Figure 10.
\endproclaim

\demo{Proof}
This is \cite{{\bf 4}, Theorem 7.1}.
\qed\enddemo

\definition{Definition 4.1.3}
A product $w_1w_2\cdots w_n$ of elements $w_i\in W(X)$ is called
{\it reduced} if $\ell(w_1w_2\cdots w_n)=\sum_i\ell(w_i)$.  We reserve
the terminology {\it reduced expression} for reduced products 
$w_1w_2\cdots w_n$ in which every $w_i \in S(X)$.

Call an element $w \in W(X)$ {\it complex} if it can be written 
as a reduced product $x_1 w_{ss'} x_2$, where $x_1, x_2 \in W(X)$ and
$w_{ss'}$ is the longest element of some rank 2 parabolic subgroup 
$\lan s, s' \ran$ such that $s$ and $s'$ correspond to adjacent nodes
in the Coxeter graph.

We define the {\it content} of $w\in W$ to be the set $c(w)$
of Coxeter generators $s\in S$ that appear in some 
(any) reduced expression for $w$.  (This can be shown not to depend 
on the reduced expression chosen, by using the theory of Coxeter groups.)

Denote by $W_c(X)$ the set of all elements of $W(X)$
that are not complex.

Let $t_w$ denote the image of the basis element $T_w \in \H(X)$ in
the quotient $TL(X)$.
\enddefinition

\proclaim{Proposition 4.1.4 (Graham)}
The set $\{ t_w : w \in W_c \}$ 
is an $\A$-basis for the algebra $TL(X)$.  \endproclaim

\demo{Proof}
See \cite{{\bf 4}, Theorem 6.2}.
\qed\enddemo

We now recall a principal result of \cite{{\bf 8}}, which establishes
the canonical basis for $TL(X)$.  This basis is a direct analogue of 
the important Kazhdan--Lusztig basis of the Hecke algebra $\H(X)$
defined in \cite{{\bf 13}}.
  
Fix a Coxeter graph, $X$.  Let 
$\A^- = \zed[v^{-1}]$, and let $\,\bar{\ }\,$ be the
involution on the ring $\A = {\Bbb Z}[v, v^{-1}]$ which satisfies 
$\bar{v} = v^{-1}$.

By \cite{{\bf 8}, Lemma 1.4}, the algebra $TL(X)$ has a $\zed$-linear 
automorphism of order $2$ that sends $v$ to $v^{-1}$ and $t_w$ to 
$t_{w^{-1}}^{-1}$.  We denote this map also by $\,\bar{\ }\,$.

Let $\L$ be the free $\A^-$-submodule of $TL(X)$ with basis
$\{\te_w : w \in W_c\}$, where $\te_w := v^{-\ell(w)} t_w$,
and let $\pi : \L \ra \L/v^{-1}\L$ be the canonical projection.

\proclaim{Proposition 4.1.5}
There exists a unique basis $\{ c_w : w \in W_c\}$ for $\L$
such that $\overline{c_w} = c_w$ and $\pi(c_w) = \pi(\te_w)$
for all $w\in W_c$.
\endproclaim

\demo{Proof}
This is \cite{{\bf 8}, Theorem 2.3}.
\qed\enddemo

The basis $\{c_w : w \in W_c\}$ is called the {\it canonical basis} (or the 
{\it IC basis}) of $TL(X)$.  It depends on the 
$t$-basis, the involution $\,\bar{\ }\,$, and the lattice $\L$.  

In this paper, we shall only be concerned with algebras $TL(X)$ of
finite rank.  We note that the canonical basis is known in many of
these cases: for types $A$, $D$ and $E$, see \cite{{\bf 8}, Theorem 3.6}; 
for types $B$ and $H$, see \cite{{\bf 9}, Theorem 2.2.1} and 
\cite{{\bf 6}, Theorem 2.1.3, Theorem 2.2.5}.  The type $I$ result follows
trivially from \cite{{\bf 9}, Proposition 1.2.3}.

It is often more convenient to work with an alternate set of algebra
generators for $TL(X)$, as follows.

\definition{Definition 4.1.6}
If $s \in S(X)$, we write $b_s \in TL(X)$ for the element $v^{-1} t_1
+ v^{-1} t_s$.  It is clear that $TL(X)$ is generated as an algebra by
the elements $b_s$.
\enddefinition

\subhead 4.2 Main results \endsubhead

It is convenient at this stage to define some named canonical basis elements of
$P_n^A$ that will turn out to be related to the elements $b_s$ of
Definition 4.1.6.

\definition{Definition 4.2.1}
Let $A$ be a normalized table algebra with basis $\BB$ and let $n, k \in \enn$.
Suppose $n > 1$ and $1 \leq k < n$.  Let $x \in \BB$.  Then the
canonical basis element $E_k(x)$ of $P_n^A$ is the one where each 
point $i$ is connected by a vertical edge to point $2n + 1 - i$,
unless $i \in \{k, k+1, 2n - k, 2n + 1 - k\}$.  Points $k$ and $k+1$
are connected by an edge, as are points $2n - k$ and $2n + 1 - k$.
All edges are labelled by $1 \in A$, except the edge connecting $k$
and $k+1$, which is labelled by $x$, and the edge connecting $2n - k$
and $2n + 1 - k$, which is labelled by $y = \bar{x}$.
\enddefinition

\example{Example 4.2.2}
Figure 11 shows the basis element $E_1(x)$ for $n = 5$.  We omit the label
$1$ (see Figure 5) since it occurs frequently.
\endexample


\topcaption{Figure 11} A basis element $E_1(x)$ \endcaption
\centerline{
\hbox to 3.027in{
\vbox to 1.375in{\vfill
        \includegraphics{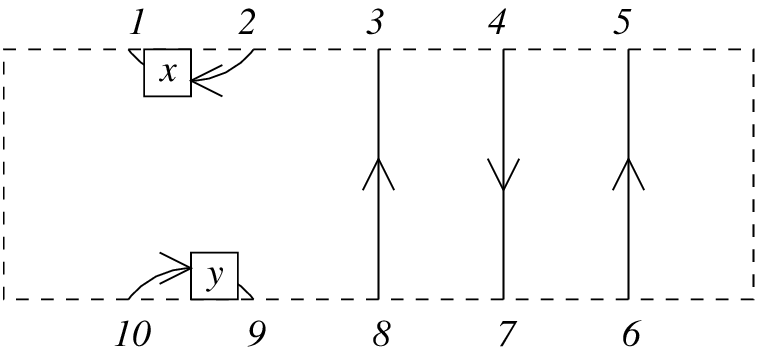}
}
\hfill}
}

To describe the sub-tabular algebras that arise from Coxeter systems
of various types, we need to define certain subsets of the canonical
basis for $P(n, r)$ (see Definition 2.4.5).  Our definitions of
$B$-admissible and $H$-admissible have previously appeared in other
forms in \cite{{\bf 6},
Definition 2.2.4} and \cite{{\bf 5}, Definition 2.2.1} respectively.

\definition{Definition 4.2.3}
We say an edge in a canonical basis element for $P(n, r)$ is {\it
$i$-decorated} if it carries a $1$-box decorated by the basis element
$u_i(x)$ of $V_r$ (see Definition 1.2.2).

\item{$\bullet$}
{A basis element, of $P(n, 3)$ is said to be {\it $B$-admissible} if it is
exposed (see Definition 2.4.3) and satisfies one of the following
three conditions:

\item{(B$1$)}{There is a $0$-decorated edge connecting points $1$ and $2n$.}
\item{(B$1'$)}{There is a $2$-decorated edge connecting points $1$ and
$2n$, and at least one non-propagating edge.}
\item{(B$2$)}{Points $1$ and $2n$ are not connected to each other and
are each the endpoint of a $1$-decorated edge, and there are no other
$1$-decorated edges.}}

\item{$\bullet$}
{A basis element, of $P(n, 4)$ is said to be {\it $H$-admissible}
if it is exposed and satisfies one of the following two conditions:

\item{(H$1$)}{If all edges are propagating, then all edges are $0$-decorated.}
\item{(H$2$)}{If not all edges are propagating, then each edge is
either $0$-decorated or $2$-decorated and both the following
statements hold:
\item{(i)}{either there is a
$2$-decorated edge connecting nodes $1$ and $2$ or there is a
$0$-decorated edge connecting nodes $i$ and $i+1$ for some $i > 1$;}
\item{(ii)}{either there is a
$2$-decorated edge connecting nodes $2n$ and $2n-1$ or there is a
$0$-decorated edge connecting nodes $2n+1-i$ and $2n-i$ for some $i > 1$.}
}}

\item{$\bullet$}
{A basis element, of $P(3, r)$ is said to be {\it $I$-admissible}
if it is exposed and satisfies one of the following three conditions:

\item{(I$1$)}{If all edges are propagating, then all edges are $0$-decorated.}
\item{(I$2$)}{A non-propagating edge is $1$-decorated if it is
transitional (see Definition 2.2.1) and $0$-decorated otherwise.}
\item{(I$3$)}{A propagating $i$-decorated edge must have $i$ odd if it
is transitional, and $i$ even if it is non-transitional.}}
\enddefinition

Note that there are a total of $2r+1$ $I$-admissible basis elements
for $P(3, r)$.

Armed with these combinatorial definitions, we can now state the main
result.

\proclaim{Theorem 4.2.4}
\item{\rm (i)}{Let $X$ be a Coxeter graph of type $A_n$.  There is a
monomorphism of $\A$-algebras $\rho_A : TL(A_n) \ra P(n+1, 2)$ such
that $\rho_A(b_i) = E_i(u_0(x))$.  The map $\rho_A$ takes canonical basis
elements of $TL(X)$ to canonical basis elements of $P(n+1, 2)$, and
its image is the Temperley--Lieb algebra, $TL(n+1, [2])$.}
\item{\rm (ii)}{Let $X$ be a Coxeter graph of type $B_n$.  There is a
monomorphism of $\A$-algebras $\rho_B : TL(B_n) \ra P(n+1, 3)$ such
that $\rho_B(b_1) = E_1(u_1(x))$ and $\rho_B(b_i) = E_i(u_0(x))$ for $i > 1$.  
The map $\rho_B$ takes canonical basis
elements of $TL(X)$ to canonical basis elements of $P(n+1, 3)$, and
its image is the subalgebra of $\D(n+1, 3)$ spanned by the
$B$-admissible basis elements.}
\item{\rm (iii)}{Let $X$ be a Coxeter graph of type $H_n$.  There is a
monomorphism of $\A$-algebras $\rho_H : TL(H_n) \ra P(n+1, 4)$ such
that $\rho_H(b_1) = E_1(u_2(x))$ and $\rho_H(b_i) = E_i(u_0(x))$ for $i > 1$.  
The map $\rho_H$ takes canonical basis
elements of $TL(X)$ to canonical basis elements of $P(n+1, 4)$, and
its image is the subalgebra of $\D(n+1, 4)$ spanned by the
$H$-admissible basis elements.}
\item{\rm (iv)}{Let $X$ be a Coxeter graph of type $I_2(m)$ for $m >
1$.  There is a
monomorphism of $\A$-algebras $\rho_I : TL(I_2(m)) \ra P(3, m-1)$ such
that $\rho_I(b_1) = E_1(u_1(x))$ and $\rho_I(b_2) = E_2(u_0(x))$.  
The map $\rho_I$ takes canonical basis
elements of $TL(X)$ to canonical basis elements of $P(3, m-1)$, and
its image is the subalgebra of $\D(3, m-1)$ spanned by the
$I$-admissible basis elements.}
\endproclaim

\demo{Proof}
We prove each part separately.

For (i), we note that the existence of $\rho_A$ and the description of
its image are well known; see \cite{{\bf 11}, Example 2.1}.  The canonical
basis of $TL(A_n)$ is shown in \cite{{\bf 8}, Theorem 3.6} to coincide
with the so-called ``monomial basis'' which in turn coincides with the
basis of diagrams of Definition 2.1.5 by \cite{{\bf 3}, Proposition
3.2.2}.  (The diagrams in \cite{{\bf 3}} are not oriented, but the
difference is cosmetic in this case.)

For (ii), we maintain the notation of Lemma 1.2.5 and 
observe using a case by case check that the
$B$-admissible basis elements span a subalgebra, $A'$, of $P(n+1, 3)$ (and
hence of $\D(n+1, 3)$, since all $B$-admissible basis elements are
exposed).  Another case by case check shows that when relation (d) of
Definition 2.2.2 is used when calculating the product of two
$B$-admissible basis elements, the label on the loop (see Figure 6) is
always a linear combination of $1$ and $z$.
By Lemma 1.2.5, there is a map, $\psi$, from $A'$ to $\D(n+1,
2)$ obtained by replacing each box decorated by $x \in V_3$ by a box
decorated by $\psi(x) \in V_2$.  Since $t(\phi(1)) = 1 = t(1')$ and
$t(\phi(z)) = 0 = t(z')$, where $t$ is the trace of Proposition 1.1.4,
it follows from Definition 2.2.2 that $\psi$ is a homomorphism.
From Definition 4.2.3, we see that the labels on a $B$-admissible 
basis element may be reconstructed
from a knowledge of which edges are not $0$-decorated.  It follows that
$\psi$ sends the canonical basis of $A'$ to a linearly independent
set in $\D(n+1, 2)$, and hence $\psi$ is injective.

The tangles of \cite{{\bf 6}, \S2.2} may be easily reconstructed from $\psi$ by
ignoring all the arrows (which convey no information as the involution
on $V_r$ is the identity map), replacing $z'$ by a square decoration
and replacing $\ds{{{1' + z'} \over 2} = {{\psi(y)} \over \sqrt{2}}}$
by a round decoration.  The relations in figures 4 and 5 of \cite{{\bf 6},
\S2.2} then correspond to the identities: $$\eqalign{
\d t(1) &= \d,\cr
\d t(\psi (y/\sqrt{2} )) &= {\d \over 2}, \cr
\psi((y/\sqrt{2})^2) &= \psi(y/\sqrt{2}),\cr
\psi(z) &= 2\psi(y/\sqrt{2}) - \psi(1),\cr
t(\psi(z)) &= 0,\cr
\psi(z)^2 &= \psi(1),\cr
\psi(y/\sqrt{2})\psi(z) &= \psi(y/\sqrt{2})
.
}$$

Using these identifications, we find that the statement of part (ii) 
is a restatement of \cite{{\bf 6}, Theorem 2.2.5} and
the map $\rho_B$ agrees with that of \cite{{\bf 6}, Theorem 2.2.3}.

The proof of (iii) is similar to, but much easier than, (ii).  The
result is a restatement of \cite{{\bf 6}, Theorem 2.1.2}, where the arrows
are ignored as before and a decorated edge in \cite{{\bf 6}, \S2.1} is
identified with a $2$-decorated edge.  The relations, shown in
\cite{{\bf 5}, Figure 5}, correspond to the relations  $$\eqalignno{
\d t(u_0(x)) &= \d, \cr
t (u_2(x)) &= 0, \cr
u_2(x)^2 &= u_0(x) + u_2(x).\cr
}$$ in $V_4$.

To prove (iv), it is convenient to construct an explicit bijection
from the set of $I$-admissible diagrams to $W_c$.  In this case, $S(I_2(m)) =
\{s_1, s_2 \}$. We send the element
described in condition (I$1$) of Definition 4.2.3 to $1 \in W_c$.
Let $D$ be a nonidentity $I$-admissible diagram.
Then $D$ has a unique propagating edge, which is $i$-decorated for
some $i$, and $D$ corresponds to the unique element $w \in W_c$ satisfying:
\item{(a)}{$s_1w < w$ if and only if points $1$ and $2$ are connected by
an edge;}
\item{(b)}{$s_2w < w$ if and only if points $2$ and $3$ are connected by
an edge;}
\item{(a$'$)}{$ws_1 < w$ if and only if points $5$ and $6$ are connected by
an edge;}
\item{(b$'$)}{$ws_2 < w$ if and only if points $4$ and $5$ are connected by
an edge;}
\item{(c)}{$\ell(w) = i+1$.}

Let $s \in S(I_2(m))$ with $\{s, s'\} = \{s_1, s_2\}$.
Using standard properties of Hecke algebras, we find that 
$c_s c_w = [2] c_w$ if $sw < w$.  Suppose now that $sw > w$.
If $\ell(w) \not\in \{m-1, 0, 1\}$, we have $$
c_s c_w = c_{sw} + c_{s'w}
.$$  If $\ell(w) = 1$ then $w = s'$ and $c_s c_{s'} = c_{s s'}$.  If
$\ell(w) = m-1$ then $sw = w_0$ and $c_s c_w = c_{s'w}$.
Using the recurrence relation in Definition
1.2.1, we see that the correspondence of the previous paragraph
defines an isomorphism of $\A$-algebras sending canonical basis
elements to $I$-admissible diagrams.
\qed\enddemo

In fact, it is possible to state a result similar to Theorem 4.2.4
more concisely, as follows.

\proclaim{Theorem 4.2.5}
Let $X$ be a (connected) Coxeter graph of rank $n > 1$ of type $A$, $B$,
$H$ or $I$, and let $m > 2$ be the highest bond label.  Then there is
a monomorphism of $\A$-algebras $\rho : TL(X) \ra P(n+1, m-1)$ such
that $\rho(b_1) = E_1(u_1(x))$ and $\rho(b_i) = E_i(u_0(x))$ for $i > 1$.  
The map $\rho$ takes canonical basis
elements of $TL(X)$ to canonical basis elements of $\D(n+1, m-1)$, and
thus $TL(X)$ equipped with its canonical basis is a sub-tabular algebra.
\endproclaim

\demo{Proof}
The assertion about sub-tabular algebras is a consequence of Theorem
3.2.3.  The other assertions, in the case of type $B$ and type $I$,
are dealt with by Theorem 4.2.4 (ii) and (iv) respectively.  

In type $A$, we take $\rho := \omega \circ \rho_A$, and in type $H$,
we take $\rho := \omega \circ \rho_H$.  By Lemma 2.4.6, $\omega$ is an
isomorphism permuting the canonical basis, so $\rho$ is a
monomorphism taking canonical basis elements to canonical basis
elements by Theorem 4.2.4.  It remains to check that $\rho$ has the
correct effect on the generators $b_i$ (in particular, on $b_1$), but 
this follows easily from the definitions of $\omega$, $\rho_A$ and $\rho_H$.
\qed\enddemo

\remark{Remark 4.2.6}
Although the statement of Theorem 4.2.5 is much shorter than that of
Theorem 4.2.4, there are situations where the results of Theorem 4.2.4
are more useful.  This is because the algebras $TL(X)$ where $X$ is of
type $A$ or $H$ are tabular algebras (see \cite{{\bf 7}, \S\S 4, 5}), and 
the table datum may be obtained easily from that of $\D(n, r)$ by 
appropriate restriction of the identifications in Theorem 4.2.4.

As well as being more uniform, the statement
of Theorem 4.2.5 gives well-defined maps for the Coxeter systems $A_2
= I_2(3)$, $B_2 = I_2(4)$ and $H_2 = I_2(5)$.

It may be easily checked that $\omega \circ \rho = \rho$ for Coxeter
systems of type $B$ or type $I_2(m)$ with $m$ even.  For $X = I_2(m)$
with $m$ odd, consideration of $\omega \circ \rho$ produces a tabular structure
for $TL(X)$ analogous to that obtained from Theorem 4.2.4 
in types $A$ and $H$.  The details are not hard to fill in.
\endremark

\subhead 4.3 Applications \endsubhead

One of the motivations for Theorem 3.2.3 relates to the following result.

\proclaim{Theorem 4.3.1 \cite{{\bf 7}}}
Let $A$ be a tabular algebra (over $\A$) with trace $\t$ and table 
datum $(\Lambda,
\Gamma, B, M, C, *)$.  Then the map $(x, y) \ra \t(xy^*)$ defines a symmetric, 
nondegenerate bilinear form on $A$ with the following properties.
\item{\rm (i)}{For all $x, y, z \in A$, $(x, yz) = (x z^*, y)$.}
\item{\rm (ii)}{The tabular basis is almost orthonormal with respect
to this bilinear form: whenever $X, X' \in \Im(C)$,
we have $$
(X, X') = 
\cases 1 \mod v^{-1} \A^- & \text{ if } X = X',\cr
0 \mod v^{-1} \A^- & \text{ otherwise.} \cr
\endcases$$}
\endproclaim

\demo{Proof}
This is \cite{{\bf 7}, Theorem 2.2.5}.
\qed\enddemo

These bilinear forms have direct relevance to generalized
Temperley--Lieb algebras, where the map $*$ is defined as follows.  

\definition{Lemma 4.3.2}
Let $X$ be an arbitrary Coxeter system.  There is a unique
$\A$-linear anti-automorphism, $*$, of $TL(X)$ that fixes the
generators $\{b_s : s \in S(X)\}$ of Definition 4.1.6.  

If $X$ is of type $A$, $B$, $H$ or $I$, the map $*$ is
induced by the tabular anti-automorphism $*$ of $P(n, r)$ via the
monomorphism $\rho$ of Theorem 4.2.5.
\enddefinition

\demo{Proof}
It is well-known that there is an $\A$-linear automorphism of $\H(X)$
sending $C'_s$ to itself for any $s \in S(X)$.  (This is the
$\A$-linear anti-automorphism that sends $T_w$ to $T_{w^{-1}}$.)
This map fixes the generators of the ideal $J(X)$ and thus gives an
automorphism of $TL(X)$ fixing $b_s = \th(C'_s)$, where $\th$ is as in
Definition 4.1.1.  Uniqueness follows as $\{b_s : s \in S(X)\}$ is a
set of algebra generators for $TL(X)$.

For the second assertion, we note that $
\rho(b_i)^* = E_i(u_{\d_{i0}})^* = E_i(u_{\d_{i0}}) = \rho(b_i)
.$
\qed\enddemo

The next result is a generalization
of the previously unproven \cite{{\bf 6}, Hypothesis 5.3.1}.

\proclaim{Corollary 4.3.3}
Let $X$ be a Coxeter system of type $A$, $B$, $H$ or $I$, and let $*$
denote the anti-automorphism of $TL(X)$ given in Lemma 4.3.2.
Then there is a 
a symmetric, nondegenerate bilinear form, $(, )$, on $TL(X)$ with the 
following properties.
\item{\rm (i)}{For all $x, y, z \in TL(X)$, $(x, yz) = (x z^*, y)$.}
\item{\rm (ii)}{The canonical basis of $TL(X)$ is almost orthonormal 
with respect to this bilinear form: whenever $w, w' \in W_c(X)$ and $c_w$ and
$c_{w'}$ are the corresponding canonical basis elements, we have $$
(c_w, c_{w'}) = 
\cases 1 \mod v^{-1} \A^- & \text{ if } w = w',\cr
0 \mod v^{-1} \A^- & \text{ otherwise.} \cr
\endcases$$}\item{\rm (iii)}
{The basis $\{\te_w : w \in W_c\}$ of $TL(X)$ defined
in \S4.1 is almost orthonormal with respect
to this bilinear form: whenever $w, w' \in W_c(X)$, we have $$
(\te_w, \te_{w'}) = 
\cases 1 \mod v^{-1} \A^- & \text{ if } w = w',\cr
0 \mod v^{-1} \A^- & \text{ otherwise.} \cr
\endcases$$}\endproclaim

\demo{Proof}
Theorem 4.2.5 shows that
the algebra $TL(X)$ is a sub-tabular algebra of a tabular algebra with
trace, with respect to the involution $*$
and its canonical basis.  The bilinear form $(, )$ is therefore inherited
from the bilinear form on $P(n, r)$ arising from theorems 3.2.3 and
4.3.1.  Theorem 4.3.1 also shows that tabular
algebras have the properties referred to in parts (i) and (ii).  The
truth of (i) and (ii) is now evident from the definition of
sub-tabular algebras.

Part (iii) follows from (ii) and Proposition 4.1.5, which shows that
$\te_w$ and $c_w$ agree modulo $v^{-1}\L$.
\qed\enddemo

Corollary 4.3.3 leads to the following characterization of the
canonical basis, up to sign.

\proclaim{Proposition 4.3.4}
Let $X$ be a Coxeter system of type $A$, $B$, $H$ or $I$, and let
$TL(X)$ be the corresponding generalized Temperley--Lieb algebra over
$\A$ equipped with the bilinear form of Corollary 4.3.3 and the
automorphism $\bar{\ }$ of \S4.1.  Suppose $x \in TL(X)$ is such that 
$\bar{x} = x$ and $(x, x) = 1 \mod v^{-1} \A^-$.  Then either
$x$ or $-x$ is a canonical basis element $c_w$.
\endproclaim

\demo{Proof}
Write $x = \sum_{w \in W_c} \l_w c_w$, where $\l_w \in \A$.  Since
$\bar{x} = x$ and $\overline{c_w} = c_w$ by Proposition 4.1.5, we must
have $\overline{\l_w} = \l_w$ for all $w$, \idest the coefficient of
$v^k$ in $\l_w$ is the same as the coefficient of $v^{-k}$.  Define
$j = j(x)$ to be $\max_{w \in W_c, \l_w \ne 0} \deg(\l_w)$; 
this is well-defined because $x \ne 0$ by the hypothesis $(x, x) = 1$.

Now $v^{-j}x = \sum_{w \in W_c} \l'_w c_w$, where $\l'_w \in \A^-$ for
all $w$.  Define $\mu_w \in \zed$ to be the constant coefficient of $\l'_w$.
Let $$
Z := \{w \in W_c: \mu_w \ne 0\} = \{w \in W_c: \l'_w \not\in v^{-1}\A^-\}
.$$  The choice of $j$ guarantees that $Z$ is nonempty.  
Since $(c_w, c_{w'}) = 0 \mod v^{-1}\A^-$ if $w \ne w'$, we have $$
(v^{-j}x, v^{-j}x) \equiv
\left(\sum_{w \in Z} (\l'_w)^2 (c_w, c_w)\right)
+
\left(\sum_{w' \in W_c \backslash Z} (\l'_{w'})^2 (c_{w'}, c_{w'})\right)
\mod v^{-1} \A^-
.$$  Applying Corollary 4.3.3 (ii), we see that
the second sum lies in $v^{-1}\A^-$ and that $$
\sum_{w \in Z} (\l'_w)^2 (c_w, c_w)
\equiv \sum_{w \in Z} (\mu_w)^2 \mod v^{-1} \A^{-1}
.$$  By hypothesis, $(x, x) = 1 \mod v^{-1} \A$, which forces $j = 0$,
$|Z| = 1$ and $\mu_w = \pm 1$ for $w \in Z$.  We conclude that $x$ is
a canonical basis element if $\mu_w = 1$, and $-x$ is a canonical
basis element if $\mu_w = -1$.
\qed\enddemo

\remark{Remark 4.3.5}
Results such as Proposition 4.3.4, where a canonical basis is
characterized up to sign by an almost orthonormality property for a
natural inner product are familiar from the work of Kashiwara 
\cite{{\bf 12}} and Lusztig \cite{{\bf 15},
Theorem 14.2.3} on quantized enveloping algebras.  
\endremark

A result similar to Proposition 4.3.4 is
true for the Kazhdan--Lusztig basis $\{C'_w: w \in W\}$ \cite{{\bf 13}} 
for the Hecke algebra, although Corollary 4.3.3 and Proposition 4.3.4
are not obvious consequences of this result.  We offer the following
conjecture concerning the role of the Kazhdan--Lusztig basis in this
context.

\proclaim{Conjecture 4.3.6}
Let $X$ be a (connected) Coxeter graph of rank $n > 1$ of type $A$, $B$,
$H$ or $I$, and let $m > 2$ be the highest bond label.  Then there is
a homomorphism of $\A$-algebras $\rho : \H(X) \ra P(n+1, m-1)$ such
that $\rho(C'_{s_1}) = E_1(u_1(x))$ and
$\rho(C'_{s_i}) = E_i(u_0(x))$ for $i > 1$.  
The map $\rho$ is injective on the set $$
\{ C'_w : w \in W(X), \ \rho(C'_w) \ne 0\}
$$ and the image of this set under $\rho$ is 
a set of canonical basis elements of $\D(n+1, m-1)$.
\endproclaim

The results of \cite{{\bf 10}, \S3.1} in conjunction with Theorem 4.2.5 
show that Conjecture 4.3.6 is true
when the Coxeter group is finite, so the only open cases are
Coxeter systems of type $H_n$ for $n > 4$.

\leftheadtext{}
\rightheadtext{}

\end
\Refs\refstyle{A}\widestnumber\key{{\bf 19}}
\leftheadtext{References}
\rightheadtext{References}

\ref\key{{\bf 1}}
\by Z. Arad and H.I. Blau
\paper On Table Algebras and Applications to Finite Group Theory
\jour J. Algebra
\vol 138 \yr 1991 \pages 137--185
\endref

\ref\key{{\bf 2}}
\by Z. Arad, E. Fisman and M. Muzychuk
\paper Generalized table algebras
\jour Isr. J. Math.
\vol 114 \yr 1999 \pages 29--60
\endref

\ref\key{{\bf 3}}
\by C.K. Fan and R.M. Green
\paper Monomials and Temperley--Lieb algebras
\jour J. Algebra
\vol 190 \yr 1997 \pages 498--517
\endref

\ref\key{{\bf 4}}
\by J.J. Graham
\book Modular representations of Hecke algebras and related algebras
\publ Ph.D. thesis
\publaddr University of Sydney
\yr 1995
\endref

\ref\key{{\bf 5}}
\by R.M. Green
\paper Cellular algebras arising from Hecke algebras of type $H_n$
\jour Math. Zeit.
\vol 229 \yr 1998 \pages 365--383
\endref

\ref\key{{\bf 6}}
\bysame
\paper Decorated tangles and canonical bases
\jour J. Algebra
\vol 246 \yr 2001 \pages 594--628
\endref

\ref\key{{\bf 7}}
\bysame
\paper Tabular algebras and their asymptotic versions
\jour J. Algebra
\miscnote in press; \hfill\newline {\tt math.QA/0107230}
\endref

\ref\key{{\bf 8}}
\by R.M. Green and J. Losonczy
\paper Canonical bases for Hecke algebra quotients
\jour Math. Res. Lett.
\vol 6 \yr 1999 \pages 213--222
\endref

\ref\key{{\bf 9}}
\bysame
\paper A projection property for Kazhdan--Lusztig bases
\jour Int. Math. Res. Not.
\vol 1 \yr 2000 \pages 23--34
\endref

\ref\key{{\bf 10}}
\bysame
\paper Fully commutative Kazhdan--Lusztig cells
\jour Ann. Inst. Fourier
\vol 51 \yr 2001 \pages 1025--1045
\endref

\ref\key{{\bf 11}}
\by V.F.R. Jones
\paper Planar Algebras, I
\miscnote preprint
\endref

\ref\key{{\bf 12}}
\by M. Kashiwara
\paper On crystal bases of the $q$-analogue of universal enveloping
algebras
\jour Duke Math. J.
\vol 63 \yr 1991 \pages 465--516
\endref

\ref\key{{\bf 13}}
\by D. Kazhdan and G. Lusztig
\paper Representations of Coxeter groups and Hecke algebras
\jour Invent. Math. 
\vol 53 \yr 1979 \pages 165--184
\endref

\ref\key{{\bf 14}}
\by M. Klazar
\paper Twelve countings with rooted plane trees
\jour European J. Combin.
\vol 18 \yr 1997
\pages 195--210
\endref

\ref\key{{\bf 15}}
\by G. Lusztig
\book Introduction to Quantum Groups
\publ Birkh\"auser \publaddr Basel \yr 1993
\endref

\ref\key{{\bf 16}}
\by P. Martin and H. Saleur
\paper The blob algebra and the periodic Temperley--Lieb algebra
\jour Lett. Math. Phys.
\vol 30 (3)
\yr 1994 
\pages 189--206
\endref

\ref\key{{\bf 17}}
\by H. Rui and C.C. Xi
\paper Cyclotomic Temperley--Lieb algebras
\miscnote preprint
\endref

\ref\key{{\bf 18}}
\by V.S. Sunder
\paper $\text{II}_1$ factors, their bimodules and hypergroups
\jour Trans. Amer. Math. Soc.
\vol 330 \yr 1992 \pages 227--256
\endref

\ref\key{{\bf 19}}
\by E. Verlinde
\paper Fusion rules and modular transformations in 2D conformal field theory
\jour Nuclear Phys. B
\vol 300 \yr 1988 \pages 360--376
\endref

\endRefs
